\documentclass{article}%
\usepackage{amssymb}
\usepackage{amsfonts}
\usepackage{graphicx}
\usepackage{amsmath}%
\providecommand{\U}[1]{\protect\rule{.1in}{.1in}}
\newtheorem{theorem}{Theorem}

\newtheorem{definition}[theorem]{Definition}

\newtheorem{lemma}[theorem]{Lemma}

\newtheorem{remark}[theorem]{Remark}

\begin{document}

\begin{center}
{\Large \textbf{Asymptotic properties for linear processes of functionals of
reversible Markov Chains}}\newline by

Magda Peligrad\footnote{\textit{ }Supported in part by a Charles Phelps Taft
Memorial Fund grant, the NSA\ grant H98230-11-1-0135} \vskip10pt
\end{center}

Department of Mathematical Sciences, University of Cincinnati, PO Box 210025,
Cincinnati, Oh 45221-0025, USA \vskip10pt

\begin{center}
\textbf{ }
\end{center}

\textbf{Abstract. }In this paper we study the asymptotic behavior of linear
processes having as innovations mean zero, square integrable functions of
stationary reversible Markov chains. In doing so we shall preserve the
generality of coefficients assuming only that they are square summable. In
this way we include in our study the long range dependence case. The only
assumption imposed on the innovations is the absolute summability of their
covariances. Besides the central limit theorem we also study the convergence
to fractional Brownian motion. The proofs are based on general results for
linear processes with stationary innovations that have interest in themselves.

\bigskip



\vskip10pt

\section{Introduction}

Let $(\xi_{i})_{i\in{\mathbb{Z}}}$ be a stationary sequence of random
variables on a probability space $(\Omega,\mathcal{K},\mathbb{P})$ with finite
second moment and zero mean $(\mathbb{E}\xi_{0}=0)$. Let $(a_{i}%
)_{i\in{\mathbb{Z}}}$ be a sequence of real numbers such that $\sum
\nolimits_{i\in{\mathbb{Z}}}a_{i}^{2}<\infty$ and denote by%
\begin{align}
X_{k}  &  =\sum_{j=-\infty}^{\infty}a_{k+j}\xi_{j}\;,\;S_{n}(X)=S_{n}%
=\sum_{k=1}^{n}X_{k},\quad\label{defX}\\
b_{n,j}  &  =a_{j+1}+\ldots+a_{j+n}\quad\mbox{and}\quad b_{n}^{2}%
=\sum_{j=-\infty}^{\infty}b_{n,j}^{2}.\nonumber
\end{align}
The linear process $(X_{k})_{k\in{\mathbb{Z}}}$ is widely used in a variety of
applied fields. It is properly defined for any square summable sequence
$(a_{i})_{i\in{\mathbb{Z}}}$ if and only if the stationary sequence of
innovations $(\xi_{i})_{i\in{\mathbb{Z}}}$ has a bounded spectral density. In
general, the covariances of $(X_{k})_{k\in{\mathbb{Z}}}$ might not be summable
so that the linear process might exhibit long range dependence.

An important theoretical question with numerous practical implications is to
prove stability of the central limit theorem under formation of linear sums.
By this we understand that if $\sum_{i=1}^{n}\xi_{i}/\sqrt{n}$ converges in
distribution to a normal variable the same holds for $S_{n}(X)\ $properly
normalized. This problem was first studied in the literature by Ibragimov
(1962) who proved that if $(\xi_{i})_{i\in{\mathbb{Z}}}$ are i.i.d. centered
with finite second moments, then $S_{n}(X)/b_{n}$ satisfies the central limit
theorem (CLT). The extra condition of finite second moment was removed by
Peligrad and Sang (2011). The central limit theorem for $S_{n}(X)/b_{n}$ for
the case when the innovations are square integrable martingale differences was
proved in Peligrad and Utev (1997) and (2006-a), where an extension to
generalized martingales was also given.

On the other hand, motivated by applications to unit root testing and to
isotonic regression, a related question is to study the limiting behavior of
$S_{[nt]}/b_{n}$ (here and throughout the paper $[x]$ denotes the integer part
of $x$). The first results for i.i.d. random innovations go back to Davydov
(1970), who established convergence to fractional Brownian motion. Extensions
to dependent settings under certain protective criteria can be found for
instance in Wu and Min (2005) and Dedecker et al. (2011), among others.

In this paper we shall address both these questions of CLT\ and convergence to
fractional Brownian motion for linear processes with functions of reversible
Markov chains innovations.

Kipnis and Varadhan (1986) considered partial sums $S_{n}$ (where $a_{0}=1,$
and $0$ elsewhere) of an additive functional zero mean of a stationary
reversible Markov chain and showed that the convergence of $var(S_{n})/n$
implies convergence of $\{S_{[nt]}/\sqrt{n},$ $0\leq t\leq1\}$ to the Brownian
motion. There is a considerable number of papers that further extend and apply
this result to infinite particle systems, random walks, processes in random
media, Metropolis-Hastings algorithms. Among others, Kipnis and Landim (1999)
considered interacting particle systems, Tierney (1994) discussed the
applications to Markov Chain Monte Carlo. Liming Wu (1999) studied the law of
the iterated logarithm.

Our first result will show that under the only assumption of absolute
summability of covariances of innovations, the partial sums of the linear
process $S_{n}(X)/b_{n}$ satisfies the central limit theorem provided
$b_{n}\rightarrow\infty.$ If we only assume the convergence of $var(S_{n})/n$
we can also treat a related linear process.

Furthermore, we shall also establish convergence to the fractional Brownian
motion under a necessary regularity condition imposed to $b_{n}^{2}.$ For a
Hurst index larger than $1/2$ we obtain a full blown invariance principle.
This is not possible without imposing additional conditions for a Hurst index
smaller than or equal to $1/2$. However we can still get the convergence of
finite dimensional distributions. For a Hurst index of $1/2$ we shall also
consider the short memory case, when the sequence of constants is absolutely
summable, and obtain convergence to the Brownian motion.

In this paper, besides a condition on the covariances, no other assumptions
such as irreducibility or aperiodicity are imposed.

The proofs are based on a result of Peligrad and Utev (2006-a) concerning the
asymptotic behavior of a class of linear processes and spectral calculus. In
addition, in Section 4.1 we develop several asymptotic results for a class of
linear processes with stationary innovations, which is not necessarily Markov
or reversible. These results have interest in themselves and can be applied to
treat other classes of linear processes. 

Applications are given to a Metropolis Hastings Markov chain, to instantaneous
functions of a Gaussian process and to random walks on compact groups.

Our paper is organized as follows: Section 2 contains the definitions, a short
background of the problem and the results. Applications are discussed in
Section 3. Section 4 is devoted to the proofs. The Appendix contains some
technical results.

\section{Definitions, background and results}

We assume that $(\gamma_{n})_{n\in\mathbb{Z}}$ is a stationary Markov chain
defined on a probability space $(\Omega,\mathcal{F},\mathbb{P})$ with values
in a general state space $(S,\mathcal{A})$. The marginal distribution is
denoted by $\pi(A)=\mathbb{P}(\gamma_{0}\in A)$. Assume that there is a
regular conditional distribution for $\gamma_{1}$ given $\gamma_{0}$ denoted
by $Q(x,A)=\mathbb{P}(\gamma_{1}\in A|\,\gamma_{0}=x)$. Let $Q$ also denotes
the Markov operator {acting via $(Qg)(x)=\int_{S}g(s)Q(x,ds)$. Next, let
$\mathbb{L}_{0}^{2}(\pi)$ be the set of measurable functions on $S$ such that
$\int g^{2}d\pi<\infty$ and $\int gd\pi=0.$ If }${g,h}\in${$\mathbb{L}_{0}%
^{2}(\pi),$ the integral }$\int_{S}g(s)h(s)d\pi$ will sometimes be denoted by
$<g,h>$.

{For some function }${g}\in${$\mathbb{L}_{0}^{2}(\pi)$, let }%
\begin{equation}
{\xi_{i}=g(\gamma}_{i}{),\ S_{n}(\xi)=\sum\limits_{i=1}^{n}\xi_{i},\ }%
\sigma_{n}({g)}=(\mathbb{E}S_{n}^{2}({\xi)})^{1/2}. \label{defcsi}%
\end{equation}
{\ Denote by $\mathcal{F}_{k}$ the $\sigma$--field generated by $\gamma_{i}$
with $i\leq k$ and by }$\mathcal{I}$ the invariant $\sigma-$field{. }

For any integrable random variable $X$ we denote $\mathbb{E}_{k}%
X=\mathbb{E}(X|\mathcal{F}_{k}).$ With this notation, $\mathbb{E}_{0}\xi
_{1}=Qg(${$\gamma$}$_{0})=\mathbb{E}(\xi_{1}|${$\gamma$}$_{0}).$ We denote by
${{||X||}_{p}}$ the norm in {$\mathbb{L}^{p}$}$(\Omega,\mathcal{F}%
,\mathbb{P}).$

The Markov chain is called reversible if $Q=Q^{\ast},$ where $Q^{\ast}$ is the
adjoint operator of $Q$. In this setting, the condition of reversibility is
equivalent to requiring that $(${$\gamma$}$_{0},${$\gamma$}$_{1})$ and
$(\gamma_{1},\gamma_{0})$ have the same distribution. Equivalently%
\[
\int_{A}Q(\omega,B)\pi(d\omega)=\int_{B}Q(\omega,A)\pi(d\omega)
\]
for all Borel sets $A,B\in\mathcal{A}$. The spectral measure of $Q$ with
respect to ${g}$ is concentrated on $[-1,1]$ and will be denoted by $\rho
_{g}.$ Then
\[
\mathbb{E}(Q^{m}g(\gamma_{0})Q^{n}g(\gamma_{0}))=<Q^{m}g,Q^{n}g>=\int
\nolimits_{-1}^{1}t^{n+m}\rho_{g}(dt).
\]
Kipnis and Varadhan (1986) assumed that%
\begin{equation}
\lim_{n\rightarrow\infty}\frac{\sigma_{n}^{2}(g)}{n}=\sigma_{g}^{2}
\label{condvar}%
\end{equation}
and proved that for any reversible ergodic Markov chain defined by
(\ref{defX}) this condition implies
\begin{equation}
W_{n}(t)=\frac{S_{[nt]}{(\xi)}}{\sqrt{n}}\Rightarrow|\sigma_{g}|W(t)\text{,}
\label{IP}%
\end{equation}
where $W(t)$ is the standard Brownian motion, $\Rightarrow$ denotes weak convergence.

As shown by Kipnis and Varadhan (1986, relation 1.1) condition (\ref{condvar})
is equivalent to
\begin{equation}
\int\nolimits_{-1}^{1}\frac{1}{1-t}\rho_{g}(dt)<\infty, \label{SR}%
\end{equation}
and then%
\[
\sigma_{g}^{2}=\int\nolimits_{-1}^{1}\frac{1+t}{1-t}\rho_{g}(dt).
\]

We shall establish the following central limit theorem:

\begin{theorem}
\label{cltlingenr} Assume that $(\xi_{j})_{j\in{\mathbb{Z}}}$ is defined by
(\ref{defcsi}) and $Q=Q^{\ast}$. Define $(X_{k})$, $S_{n}$ and $b_{n}$ as in
(\ref{defX}). Assume that $b_{n}\rightarrow\infty$ as $n\rightarrow\infty$
and
\begin{equation}
\sum_{j\geq0}|\mathrm{cov}(\xi_{0},\xi_{j})|<\infty.\label{abscov}%
\end{equation}
Then, there is a nonnegative random variable $\eta$ measurable with respect to
$\mathcal{I}$ such that $n^{-1}\mathbb{E}((\sum_{k=1}^{n}\xi_{k}%
)^{2}|\mathcal{F}_{0})\rightarrow\eta$\ in $L_{1}$ as $n\rightarrow\infty$ and
$\mathbb{E}\eta=\sigma_{g}^{2}.$ In addition
\[
\lim_{n\rightarrow\infty}\frac{\mathrm{Var}(S_{n}(X))}{b_{n}^{2}}=\sigma
_{g}^{2}%
\]
and
\begin{equation}
\frac{S_{n}(X)}{b_{n}}\Rightarrow\sqrt{\eta\ }N\text{ as }n\rightarrow
\infty,\label{CLT}%
\end{equation}
where $N$ is a standard normal variable independent on $\eta.$ Moreover if the
sequence $(\xi_{i})_{i\in{\mathbb{Z}}}$ is ergodic the central limit theorem
in (\ref{CLT}) holds with $\eta=\sigma_{g}^{2}.$
\end{theorem}

It should be noted that under the conditions of this theorem $\sigma_{g}^{2}$
also has the following interpretation: the stationary sequence $(\xi
_{i})_{i\in{\mathbb{Z}}}$ has a continuous spectral density $f(x)$ and
$\sigma_{g}^{2}=2\pi f(0).$

\bigskip

In order to present the functional form of the CLT\ we introduce a regularity
assumption which is necessary for this type of result. We denote by $D([0,1])$
the space of functions defined on $[0,1]$ which are right continuos and have
left hand limits at any point.

\begin{definition}
We say that a positive sequence $(b_{n}^{2})_{n\geq1}$ is regularly varying
with exponent $\beta>0$ if for any $t\in]0,1]$,
\begin{equation}
\frac{b_{[nt]}^{2}}{b_{n}^{2}}\rightarrow t^{\beta}\,\text{{as }%
}\,n\rightarrow\infty. \label{hyposn}%
\end{equation}

\end{definition}

We shall separate the case $\beta\in]1,2]$ from the case $\beta\in]0,1].$

\begin{theorem}
\label{IP1}Assume that the conditions of Theorem \ref{cltlingenr} are
satisfied and in addition $b_{n}^{2}$, defined by (\ref{defX}), is regularly
varying with exponent $\beta$ for a certain $\beta\in]1,2]$. Then, the process
$\{b_{n}^{-1}S_{[nt]}(X),t\in\lbrack0,1]\}$ converges in $D([0,1])$ to
$\sqrt{\eta}W_{H}$ where $W_{H}$ is a standard fractional Brownian motion
independent of $\eta$ with Hurst index $H=\beta/2$.
\end{theorem}

The case $\beta\in]0,1]$ is more delicate. For this case we only give the
convergence of the finite dimensional distributions since there are
counterexamples showing that the tightness might not hold without additional
assumptions. As a matter of fact, for $\beta=1$, it is known from
counterexamples given in Wu and Woodroofe (2004) and also in Merlev\`{e}de and
Peligrad (2006) that the weak invariance principle may not be true for the
partial sums of the linear process with i.i.d. square integrable innovations.

\begin{theorem}
\label{IP2} Assume that the conditions of Theorem \ref{cltlingenr} are
satisfied and in addition $b_{n}^{2}$ is regularly varying with exponent
$\beta$ for a certain $\beta\in]0,1]$. Then the finite dimensional
distributions of $\{b_{n}^{-1}S_{[nt]},t\in\lbrack0,1]\}$ converges to the
corresponding ones of $\sqrt{\eta}W_{H}$, where $W_{H}$ is a standard
fractional Brownian motion independent of $\eta$ with Hurst index $H=\beta/2$.
\end{theorem}

In the context of Theorems \ref{IP1} and \ref{IP2}, condition (\ref{hyposn})
is necessary for the conclusion of this theorem (see Lamperti, 1962). This
condition has been also imposed by Davydov (1970) for studying the weak
invariance principle of linear processes with i.i.d. innovations.

The following theorem is obtained under condition (\ref{condvar}). 

\begin{theorem}
\label{blocks-inn}Assume that $(\xi_{j})$ is defined by \ref{defcsi}\ and
condition (\ref{condvar}) is satisfied. Define%
\begin{equation}
X_{k}^{\prime}=\sum_{j=-\infty}^{\infty}a_{k+j}(\xi_{j}+\xi_{j+1}%
)\;,\;S_{n}(X^{\prime})=\sum_{k=1}^{n}X_{k}^{\prime},.\label{defX'}%
\end{equation}
Then the conclusion of Theorems \ref{cltlingenr}, \ref{IP1} and \ref{IP2} hold
for $S_{n}(X^{\prime})$. In this case $\eta$ is identified as the limit
$n^{-1}\mathbb{E}(\sum_{k=1}^{n}(\xi_{k}+\xi_{k+1})^{2}|\mathcal{F}%
_{0})\rightarrow\eta$\ in $L_{1}$ as $n\rightarrow\infty$. Furthermore, the
stationary sequence $(\xi_{k}+\xi_{k+1})_{k\in{\mathbb{Z}}}$ has a continuous
spectral density $h(x)$ and $\mathbb{E}\eta=2\pi h(0)=\lim_{n\rightarrow
\infty}\mathrm{Var}S_{n}(X^{\prime})/b_{n}^{2}.$
\end{theorem}

\bigskip

We shall present next the short memory case:

\begin{theorem}
\label{shortmem} Assume now $\sum\nolimits_{i\in{\mathbb{Z}}}|a_{i}|<\infty$
and let $(X_{k})_{k\geq1}$ be as in Theorem \ref{cltlingenr}. Assume that
condition (\ref{condvar}) is satisfied. Then the process $\{S_{[nt]}/\sqrt
{n},t\in\lbrack0,1]\}$ converges in $D([0,1])$ to $\sqrt{\eta}|A|W$ where $W$
is a standard Brownian motion and $A=\sum\nolimits_{i\in{\mathbb{Z}}}a_{i}$.
\end{theorem}

\begin{remark}
\label{Block}It is easy to see that Theorems (\ref{shortmem}) extends Kipnis
Varadhan result to linear processes. (\ref{condvar}) 
\end{remark}

We give a few examples of sequences $(a_{n})$ satisfying the conditions of our
theorems. In these examples the notation $a_{n}\sim b_{n}$ means $a_{n}%
/b_{n}\rightarrow1$ as $n\rightarrow\infty.$

\medskip

\noindent\textbf{Example 1}. For the selection $a_{i}\sim i^{-\alpha}\ell(i)$
where $\ell$ is a slowly varying function at infinity and $1/2<\alpha<1$ for
$i\geq1$ and $a_{i}=0$ elsewhere, then, $b_{n}^{2}\sim\kappa_{\alpha
}n^{3-2\alpha}\ell^{2}(n)$ (see for instance Relations (12) in Wang \textit{et
al.} (2003)), where $\kappa_{\alpha}$ is a positive constant depending on
$\alpha$. Clearly, Theorem \ref{IP1} applies.

\medskip

\noindent\textbf{Example 2}. Let us consider now the fractionally integrated
processes since they play an important role in financial time series modeling
and they are widely studied. Such processes are defined for $0<d<1/2$ by
\begin{equation}
X_{k}=(1-B)^{-d}\xi_{k}=\sum_{i\geq0}a_{i}\xi_{k-i}\ \mbox{ with }a_{i}%
=\frac{\Gamma(i+d)}{\Gamma(d)\Gamma(i+1)}\,, \label{deffractlin}%
\end{equation}
where $B$ is the backward shift operator, $B\varepsilon_{k}=\varepsilon_{k-1}%
$. For this example, by the well known fact that for any real $x,$
$\lim_{n\rightarrow\infty}\Gamma(n+x)/n^{x}\Gamma(n)=1,$ we have$\ \lim
_{n\rightarrow\infty}a_{n}/n^{d-1}=1/\Gamma(d)$. Theorem \ref{IP1} applies
with $\beta=2d+1$, since for $k\geq1$ we have $a_{k}\sim\kappa_{d}k^{d-1}$ for
some $\kappa_{d}>0$ and $a_{k}=0$ elsewhere.

\medskip

\noindent\textbf{Example 3}. Now, if we consider the following selection of
$(a_{k})_{k\geq0}$: $a_{0}=1$ and $a_{i}=(i+1)^{-\alpha}-i^{-\alpha}$ for
$i\geq1$ with $\alpha\in]0,1/2[$ and $a_{i}=0$ elsewhere, then Theorem
\ref{IP2} applies. Indeed for this selection, $b_{n}^{2}\sim\kappa_{\alpha
}n^{1-2\alpha}$, where $\kappa_{\alpha}$ is a positive constant depending on
$\alpha$.

\medskip

\noindent\textbf{Example 4}. Finally, if $a_{i}\sim i^{-1/2}(\log i)^{-\alpha
}$ for some $\alpha>1/2$, then $b_{n}^{2}\sim n^{2}(\log n)^{1-2\alpha
}/(2\alpha-1)$ (see Relations (12) in Wang \textit{et al.} (2003)). Hence
(\ref{hyposn}) is satisfied with $\beta=2$.

\section{Applications}

\subsection{Application to a Metropolis Hastings Markov chain.}

In this subsection we analyze a standardized example of a stationary
irreducible and aperiodic Metropolis-Hastings algorithm with uniform marginal
distribution. This type of Markov chain is interesting since it can easily be
transformed into Markov chains with different marginal distributions. Markov
chains of this type are often studied in the literature from different points
of view. See, for instance, Doukhan et al (1994) and Longla et al (2012) among
many others.

Let $E=[-1,1]$ and let $\upsilon$ be a symmetric atomless law on $E$. The
transition probabilities are defined by
\[
Q(x,A)=(1-|x|)\delta_{x}(A)+|x|\upsilon(A),
\]
where $\delta_{x}$ denotes the Dirac measure. Assume that $\theta=\int
_{E}|x|^{-1}\upsilon(dx)<\infty$. Then there is a unique invariant measure
\[
\pi(dx)=\theta^{-1}|x|^{-1}\upsilon\,(dx)
\]
and the stationary Markov chain $(\gamma_{k})$ generated by $Q(x,A)$ and $\pi$
is reversible and positively recurrent, therefore ergodic.

\begin{theorem}
Let $g(-x)=-g(x)$ for any $x\in E$ and assume
\[
\int_{0}^{1}g^{2}(x)x^{-2}dv<\infty.
\]
Then, the conclusions of all our theorems in Section 2 hold for $(X_{k})$ and
$S_{n}(X)$ defined by (\ref{defX}) with
\[
\eta=\sigma_{g}^{2}=\theta^{-1}(\int_{E}g^{2}(x)|x|^{-1}\upsilon
(dx)\,+2\int_{E}g^{2}(x)|x|^{-2}\upsilon(dx)).
\]

\end{theorem}

\textbf{Proof}. Since $g$ is an odd function we have
\begin{equation}
\mathbb{E}(g(\gamma_{k})|\gamma_{0})=(1-|\gamma_{0}|)^{k}g(\gamma_{0})\text{
a.s.} \label{operator}%
\end{equation}
Therefore, for any $j\geq0$,
\[
\mathbb{E}(X_{0}X_{j})=\mathbb{E}(g(\gamma_{0})\mathbb{E}(g(\gamma_{j}%
)|\gamma_{0}))=\theta^{-1}\int_{E}g^{2}(x)(1-|x|)^{j}|x|^{-1}\upsilon(dx).
\]
Then,
\begin{equation}
\sum_{j=1}^{k-1}|\mathbb{E}(X_{0}X_{j})|\leq2\theta^{-1}\sum_{j=1}^{k-1}%
\int_{0}^{1}g^{2}(x)(1-x)^{j}x^{-1}\upsilon(dx)\leq2\theta^{-1}\int_{0}%
^{1}g^{2}(x)x^{-2}\upsilon(dx) \label{estimatevar}%
\end{equation}
and therefore condition (\ref{abscov}) is satisfied. $\ \Diamond$

\subsection{\textbf{Linear process of instantaneous functions of a Gaussian
sequence}}

\begin{theorem}
Let $({\xi}_{k})_{k\in{\mathbb{Z}}}$ be instantaneous functions of a
stationary Markov Gaussian sequence $(\gamma_{n}),$ ${\xi}_{k}=g(\gamma_{n})$
where $g$ is a measurable real function such that $\mathbb{E}g(\gamma_{n})=0$
and $\mathbb{E}g^{2}(\gamma_{n})<\infty.$ Define $X_{k}$ and $S_{n}(X)$ by
(\ref{defX}). Then the conclusion of our theorems in Section 2 hold.
\end{theorem}

\textbf{Proof}. In order to apply our results, because $(\gamma_{n})$ is
reversible, we have only to check condition (\ref{abscov}). Under our
conditions $g$ can be expanded in Hermite polynomials $g(x)=\sum_{j\geq1}%
c_{j}H_{j}(x)$, where $\sum_{j=1}c_{j}^{2}j!<\infty.$

For computing the covariances we shall apply the following well-known formula:
if $a$ and $b$ are jointly Gaussian random variables, $\mathbb{E}%
a=\mathbb{E}b=0$, $\mathbb{E}a^{2}=\mathbb{E}b^{2}=1$, $r=\mathbb{E}ab$, then%

\[
\mathbb{E}H_{k}(a)H_{l}(b)=\delta(k,l)r^{k}k!\text{,}%
\]
where $\delta$ denotes the Kronecker delta. It follows that
\[
cov({\xi}_{0},{\xi}_{k})=\mathbb{E}\sum_{j\geq1}c_{j}^{2}H_{j}(\gamma
_{0})H_{j}(\gamma_{k})=\sum_{j\geq1}c_{j}^{2}r_{k}^{j}j!.
\]
Clearly, because under our condition it is known that $r_{k}=\exp(-\alpha
k/2)$ for some $\alpha>0,$ then
\[
|cov({\xi}_{0},{\xi}_{k})|\leq\exp(-\alpha k/2)\sum_{j\geq1}c_{j}^{2}j!
\]
and the result follows. $\ \Diamond$

\bigskip

For a particular class of weights of the form in Example 3, we mention that
Breuer and Major (1983) studied this problem\ for Gaussian chains without
Markov assumption.

\subsection{Application to random walks on compact groups}

In this section we shall apply our results to random walks on compact groups.

Let $\mathcal{X}$ be a compact abelian group, $\mathcal{A}$ a sigma algebra of
Borel subsets of $\mathcal{X}$ and $\pi$ the normalized Haar measure on
$\mathcal{X}$. The group operation is denoted by $+$. Let $\nu$ be a
probability measure on $(\mathcal{X},\mathcal{A)}$. The random walk on
$\mathcal{X}$ defined by $\nu$ is the stationary Markov chain having the
transition function%
\[
(x,A)\rightarrow Q(x,A)=\nu(A-x)\text{.}%
\]
The corresponding Markov operator denoted by $Q$ is defined by%
\[
(Qf)(x)=f\ast\nu(x)=\int_{\mathcal{X}}f(x+y)\nu(dy)\text{.}%
\]
The Haar measure is invariant under $Q.$ We shall assume that $\nu$ is not
supported by a proper closed subgroup of $\mathcal{X},$ a condition that is
equivalent to $Q$ being ergodic. In this context%
\[
(Q^{\ast}f)(x)=f\ast\nu^{\ast}(x)=\int_{\mathcal{X}}f(x-y)\nu(dy)\text{,}%
\]
where $\nu^{\ast}$ is the image of measure $\nu$ by the map $x\rightarrow-x.$
Thus $Q$ is symmetric on $\mathbb{L}_{2}(\pi)$ if and only if $\nu$ is
symmetric on $\mathcal{X}$, that is $\nu=\nu^{\ast}.$

The dual group of $\mathcal{X}$, denoted by $\mathcal{\hat{X}}$, is discrete.
Denote by $\hat{\nu}$ the Fourier transform of the measure $\nu,$ that is the
function%
\[
g\rightarrow\hat{\nu}(g)=\int_{\mathcal{X}}g(x)\nu(dx)\text{ }\ \text{with
}g\in\mathcal{\hat{X}}\text{.}%
\]
A function $f\in${$\mathbb{L}^{2}(\pi)$} has the Fourier expansion%
\[
f=\sum\limits_{g\in\mathcal{\hat{X}}}\hat{f}(g)g\text{.}%
\]
Ergodicity of $Q$ is equivalent to $\hat{\nu}(g)\neq1$ for any non-identity
$g\in\mathcal{\hat{X}}.$ By arguments in Borodin and Ibragimov (1994, ch. 4,
section 9) and also Derriennic and Lin (2001, Section 8) condition (\ref{SR})
takes the form%
\begin{equation}
\sum_{1\neq g\in\mathcal{\hat{X}}}\frac{|\hat{f}(g)|^{2}}{|1-\hat{\nu}%
(g)|}<\infty\text{.} \label{G1}%
\end{equation}

Combining these considerations with the results in Section 2 we obtain the
following result:

\begin{theorem}
Let $\nu$ be ergodic and symmetric on $\mathcal{X}$. Let $(\xi_{i})$ be the
stationary Markov chain with marginal distribution $\pi$ and transition
operator $Q$. If for $g$ in {$\mathbb{L}_{0}^{2}(\pi)$ condition }(\ref{G1})
is satisfied then the conclusions of Theorem \ref{blocks-inn} in Section 2
hold for $(X_{k}^{\prime})$ and $S_{n}(X^{\prime})$ defined by (\ref{defX'}).
\end{theorem}

\section{Proofs}

\subsection{Preliminary general results\label{prem}}

This section contains some general results for linear processes of stationary
sequences which are not necessarily Markov. We start by mentioning the
following theorem which is a variant of a result from Peligrad and Utev
(2006-a). See also Proposition 5.1 in Dedecker et al. (2011).

\begin{theorem}
\label{cltlingen}Let $(\xi_{k})_{k\in\mathbb{Z}}$ be a strictly stationary
sequence of centered square integrable random variables such that%
\begin{equation}
\Gamma_{j}=\sum_{k=0}^{\infty}|\mathbb{E(}\xi_{j+k}\mathbb{E}_{0}\xi
_{j})|<\infty\text{ and}\;\frac{1}{p}\sum_{j=1}^{p}\Gamma_{j}\rightarrow
0\text{ as}\;p\rightarrow\infty. \label{Mgen}%
\end{equation}
For any positive integer $n$, let $(d_{n,i})_{i\in\mathbb{Z}}$ be a triangular
array of numbers satisfying, for some positive $c$,
\begin{equation}
\sum_{i\in\mathbb{Z}}d_{n,i}^{2}\rightarrow c^{2}\,\text{and}\,\sum
_{j\in\mathbb{Z}}^{\ }(d_{n,j}-d_{n,j-1})^{2}\rightarrow0\,\text{as
}n\rightarrow\infty. \label{A}%
\end{equation}
In addition assume
\begin{equation}
\sup_{j\in\mathbb{Z}}|d_{n,j}|\rightarrow0\,\ {\text{as}}\text{ \ }%
n\rightarrow\infty. \label{B}%
\end{equation}
Then $\sum_{j\in\mathbb{Z}}d_{n,j}\xi_{j}$ converges in distribution to
$\sqrt{\eta}cN$ where $N$ is a standard Gaussian random variable independent
of $\eta$. The variable $\eta$ is measurable with respect to the invariant
sigma field $\mathcal{I}$ and $n^{-1}\mathbb{E}((\sum_{k=1}^{n}\xi_{k}%
)^{2}|\mathcal{F}_{0})\rightarrow\eta\quad$\ in $L_{1}$ as$\quad
n\rightarrow\infty$. Furthermore $(\xi_{i})_{i\in{\mathbb{Z}}}$ has a
continuous spectral density $f(x)$ and $\ \mathbb{E}\eta=2\pi f(0).$ If the
sequence $(\xi_{i})_{i\in{\mathbb{Z}}}$ is ergodic we have $\eta=2\pi f(0).$
\end{theorem}

\textbf{Proof}. The proof follows the lines of Theorem 1 from Peligrad and
Utev (2006-a). We just have to repeat the arguments there with $b_{n,i}/b_{n}$
replaced by $d_{n,i}$ and take into account that the properties (\ref{A}) and
(\ref{B}) are precisely all is needed to complete the proof. $\ \Diamond$

\bigskip

Next we shall establish the convergence of finite dimensional distributions.

\begin{theorem}
\label{IPgen}Define $(X_{k})$ and $S_{n}$ by (\ref{defX}) and assume condition
(\ref{Mgen}) is satisfied. Then $S_{n}/b_{n}$ converges in distribution to
$\sqrt{\eta}N$ where $N$ and $\eta$ are as in Theorem \ref{cltlingen}. If we
assume in addition that condition (\ref{hyposn}) is satisfied, then the finite
dimensional distributions of $\{W_{n}(t)=b_{n}^{-1}S_{[nt]},t\in\lbrack0,1]\}$
converge to the corresponding ones of $\sqrt{\eta}W_{H}$, where $W_{H}$ is a
standard fractional Brownian motion independent of $\eta$ with Hurst index
$H=\beta/2.$
\end{theorem}

\textbf{Proof}. The central limit theorem part requires just to verify the
conditions of Theorem \ref{cltlingen} for $d_{n,j}=b_{n,j}/b_{n}$ and $c=1.$
Condition (\ref{B}) was verified in Peligrad and Utev (1997, page 448-449)
while condition (\ref{A}) was verified in Lemma A.1. in Peligrad and Utev (2006-a).

We shall prove next the second part of the theorem. Notice that if we impose
(\ref{hyposn}), for each $t$ fixed \thinspace%
\begin{equation}
\text{var}(W_{n}(t))\rightarrow2\pi f(0)t^{\beta} \label{conv1}%
\end{equation}
and $W_{n}(t)\Rightarrow\eta t^{\beta/2}N.$

Let $0\leq t_{1}\leq...\leq t_{k}\leq1.$ By Cram\`{e}r-Wold device, in order
to find the limiting distribution of $(W_{n}(t_{i}))_{1\leq i\leq k}$ we have
to study $V_{n}=\sum_{i=1}^{k}u_{i}W_{n}(t_{i})$ where $u_{i}$ is a real
vector$.$ Let us compute its limiting variance. To find it, let $0\leq s\leq
t\leq1.$ By using the fact that for any two real numbers $a$ and $b$ we have
$a(a-b)=(a^{2}+(a-b)^{2}-b^{2})/2,$ we obtain the representation:%
\begin{gather*}
cov(W_{n}(t),W_{n}(s))=var(W_{n}(s))+cov(W_{n}(s),W_{n}(t)-W_{n}(s))\\
=var(W_{n}(s))+1/2[var(W_{n}(t)-W_{n}(s))+var(W_{n}(t))-var(W_{n}(s))].
\end{gather*}
By stationarity,%
\[
var(W_{n}(t)-W_{n}(s))=var(W_{[nt]-[ns]}),
\]
and by (\ref{conv1}) and the fact that $b_{n}\rightarrow\infty$ we obtain%
\begin{equation}
\lim_{n\rightarrow\infty}cov(W_{n}(t),W_{n}(s))=\pi f(0)(s^{\beta}+t^{\beta
}-|t-s|^{\beta}). \label{conv2}%
\end{equation}
So,%
\begin{equation}
\lim_{n\rightarrow\infty}\frac{1}{2\pi f(0)}\text{var}(V_{n})=\sum_{i=1}%
^{k}u_{i}^{2}t_{i}^{\beta}+\sum_{i=1}^{k-1}\sum_{j=i+1}^{k}u_{i}u_{j}%
(t_{i}^{\beta}+t_{j}^{\beta}-(t_{j}-t_{i})^{\beta})=B_{k}. \label{var1}%
\end{equation}
Writing now
\[
V_{n}=\sum_{i=1}^{k}u_{i}W_{n}(t_{i})=\sum_{j\in\mathbb{Z}}d_{n,j}(k)\xi_{j},
\]
where $d_{n,j}(k)=\sum_{i=1}^{k}u_{i}b_{[nt_{i}],j}/b_{n},$ we shall apply
Theorem \ref{cltlingen}. The second part of (\ref{A}) and (\ref{B}) were
verified in Peligrad and Utev (1996 and 2006-a). It remains to verify the
first part of condition (\ref{A}). By the point (iii) of Lemma \ref{spec} in
the Appendix we obtain%
\[
\text{var}(V_{n})/\sum_{j\in\mathbb{Z}}d_{n,j}^{2}(k)\rightarrow2\pi f(0),
\]
which combined with (\ref{var1}) implies that the first part of (\ref{A}) is
verified with $c^{2}=\lim_{n\rightarrow\infty}\sum_{j\in\mathbb{Z}}d_{n,j}%
^{2}(k)=B_{k}.$ In other words, the finite dimensional distributions are
convergent to those of a fractional Brownian motion with Hurst index
$\beta/2.$ $\ \Diamond$

\bigskip

\textbf{Discussion on tightness. }As we mentioned above,\textbf{ }for
$\beta\leq1$ the conditions of Theorem \ref{IPgen} are not sufficient to imply tightness.

However for $\beta>1$ we can obtain tightness in $D([0,1])$ endowed with
Skorohod topology. By the point (i) of Lemma \ref{spec} in Appendix we have
the inequality
\[
{\mathbb{E}}|S_{k}|^{2}\leq\left(  \mathbb{E}[\xi_{0}^{2}]+2\sum
_{k\in\mathbb{Z}}|\mathbb{E}(\xi_{0}\xi_{k})|\right)  \sum_{j\in\mathbb{Z}%
}b_{k,j}^{2}.
\]
Therefore, by using (\ref{Mgen}) and (\ref{hyposn}), the conditions of Lemma
2.1 p. 290 in Taqqu (1975) are satisfied when $\beta>1$, and the tightness
follows. $\ \Diamond$

\bigskip

To treat the short memory case we mention the following result in Peligrad and
Utev (2006-b).

\begin{theorem}
\label{short memory}Assume that $X_{k\text{ }}$and $S_{n}$ are defined by
(\ref{defX}) and $\sum_{i\in Z}|a_{i}|<\infty$. Moreover assume that for some
$c_{n}>0$ the innovations satisfy the invariance principle%
\[
c_{n}^{-1}S_{[nt]}(\xi)\Rightarrow\eta W(t),
\]
where $\eta$ is $\mathcal{I}$--measurable and $W$ is a standard Brownian
motion on $[0,1]$ independent on $\mathcal{I}.$ In addition assume that the
following condition holds:
\begin{equation}
{\mathbb{E}}\max_{1\leq j\leq n}|S_{j}(\xi)|\leq Cc_{n}. \label{key1}%
\end{equation}
where $C$ is a positive constant. Then, the linear process also satisfies the
invariance principle, i.e. $c_{n}^{-1}S_{[nt]}(X)\Rightarrow\eta|A|W(t)$ as
$n\rightarrow\infty$ where $A=\sum_{i\in\mathbb{Z}}a_{i}$.
\end{theorem}

\subsection{ Normal and reversible Markov Chains}

In this subsection we give the proofs of the theorems stated in Section 2. The
goal is to verify condition (\ref{Mgen}) that will assure that all the results
in the subsection \ref{prem} are valid.

We start by applying the general results to normal Markov chains, for which
$QQ^{\ast}=Q^{\ast}Q.$ For this case condition (\ref{Mgen}) is implied by%
\begin{equation}
\sum_{k\geq0}||Q^{k}g||_{2}^{2}<\infty\text{.} \label{cond normal}%
\end{equation}
Indeed, we start by rewriting (\ref{Mgen}) in operator notation:%
\begin{gather*}
|\mathbb{E}[\xi_{j+k}\mathbb{E}(\xi_{j}|\mathcal{F}_{0})]|=|\mathbb{E}%
(\mathbb{E}_{0}\xi_{k+j}\mathbb{E}_{0}\xi_{j})|=|<Q^{k+j}g,Q^{j}g>|=\\
|<Q^{[k/2]+j}g,(Q^{\ast})^{k-[k/2]}Q^{j}g>|\leq||Q^{[k/2]+j}g||_{2}||(Q^{\ast
})^{k-[k/2]}Q^{j}g||_{2}.
\end{gather*}
For normal operator, by using the properties of conditional expectation, we
have%
\[
||(Q^{\ast})^{k-[k/2]}Q^{j}g||_{2}=||Q^{j}(Q^{\ast})^{k-[k/2]}g||_{2}%
\leq||(Q^{\ast})^{k-[k/2]}g||_{2}.
\]
Since for all $\varepsilon>0,$ and any two numbers $a$ and $b$ we have
$|ab|\leq a^{2}/2\varepsilon+\varepsilon b^{2}/2$, by the above considerations
we easily obtain
\begin{align*}
\sum_{k\geq0}|\mathbb{E}[\xi_{j+k}\mathbb{E}(\xi_{j}|\mathcal{F}_{0})]|  &
\leq\sum_{k\geq0}||Q^{[k/2]+j}g||_{2}||Q^{k-[k/2]}g||_{2}\\
&  \leq\frac{1}{\varepsilon}\sum_{k\geq j}||Q^{k}g||_{2}^{2}+\varepsilon
\sum_{k\geq0}||Q^{k}g||_{2}^{2},
\end{align*}
condition (\ref{Mgen}) is verified under (\ref{cond normal}), by letting
$j\rightarrow\infty$ followed by $\varepsilon\rightarrow0.$

In terms of spectral measure $\rho_{g}(dz),$ condition (\ref{cond normal}) is
implied by%
\[
\int_{D}\frac{1}{1-|z|}\rho_{g}(dz)<\infty,
\]
where $D$ is the unit disk. Note that this condition is stronger than the
condition needed for the validity of CLT\ for the partial sums (i.e. the case
$a_{1}=1,a_{i}=0$ elsewhere), which requires only the condition $\int_{D}%
\frac{1}{|1-z|}\rho_{g}(dz)<\infty$ (see Gordin and Lifshitz (1981), or in Ch.
IV in \ Borodin and Ibragimov (1994)).

For the reversible Markov chains just notice that
\[
\mathbb{E}[\xi_{j+k}\mathbb{E}(\xi_{j}|\mathcal{F}_{0})]=\int_{-1}^{1}%
t^{2j+k}\rho_{g}(dz)=cov(\xi_{0},\xi_{2j+k})
\]
and then, condition (\ref{Mgen}) is verified under (\ref{abscov}) because%
\[
\sum_{k\geq0}|\mathbb{E}[\xi_{j+k}\mathbb{E}(\xi_{j}|\mathcal{F}_{0}%
)]|=\sum_{k\geq2j}|cov(\xi_{0},\xi_{k})|\rightarrow0\text{ as }j\rightarrow
\infty.
\]
Theorems \ref{cltlingenr}, \ref{IP1} and \ref{IP2} follow as simple
applications of the results in subsection \ref{prem}.

\bigskip

\textbf{Proof of Theorem \ref{blocks-inn} }

In order to prove this theorem, we shall also apply Theorem \ref{IPgen} along
to the tightness discussion at the end of Section 3. We denote $\gamma_{j}%
=\xi_{j}+\xi_{j+1}$ and verify condition (\ref{Mgen}) for this sequence of
innovations. We have
\[
|\mathbb{E}(\gamma_{k+j}\mathbb{E}_{0}\gamma_{j})|=|<Q^{k+j}g+Q^{k+j+1}%
g,Q^{j}g+Q^{j+1}g>|
\]
and by spectral calculus%
\[
\sum_{k\geq0}|<Q^{k+j}g+Q^{k+j+1}g,Q^{j}g+Q^{j+1}g>|=\sum_{k\geq0}|\int
_{-1}^{1}t^{k+2j}(1+t)^{2}d\rho_{g}|.
\]
We divide the sum in $2$ parts, according to $k$ even or odd. When $k=2u$ the
sum has positive terms and it can be written as
\[
\sum_{u\geq0}\int_{-1}^{1}t^{2u+2j}(1+t)^{2}d\rho_{g}\leq\int_{-1}^{1}%
\frac{t^{2j}}{1-t^{2}}(1+t)^{2}d\rho_{g}=\int_{-1}^{1}\frac{t^{2j}(1+t)}%
{1-t}d\rho_{g}.
\]
When $k$ is odd
\begin{gather*}
\sum_{k\geq1,k\text{ odd}}|\int_{-1}^{1}t^{k+2j}(1+t)^{2}d\rho_{g}|\leq
\int_{-1}^{1}\sum_{k\geq1,k\text{ odd}}|t^{k+2j}(1+t)^{2}|d\rho_{g}\\
\leq\int_{-1}^{1}\sum_{k\geq1,k\text{ odd}}|t^{k-1+2j}(1+t)^{2}|d\rho_{g}%
\leq\sum_{u\geq0}|t^{2u+2j}(1+t)^{2}|d\rho_{g},
\end{gather*}
and we continue the computation as for the case $k$ even. It follows that
\[
\frac{1}{m}\sum_{j=1}^{m}\sum_{k\geq0}|\mathbb{E}(\gamma_{k+j}\mathbb{E}%
_{0}\gamma_{j})|\leq\frac{2}{m}\sum_{j=1}^{m}\int_{-1}^{1}\frac{t^{2j}%
(1+t)}{1-t}d\rho_{g}.
\]
Note that (\ref{SR}) implies that $\rho_{g}(1)=0.$ We also have $m^{-1}%
\sum_{j=1}^{m}t^{2j}(1+t)$ is convergent to $0$ for all $t\in\lbrack-1,1)$.
Furthermore, $m^{-1}\sum_{j=1}^{m}t^{2j}(1+t)$ is dominated by $2$ and in view
of (\ref{SR}) and Lebesgue dominated convergence theorem we have
\[
\lim_{m\rightarrow\infty}\int_{-1}^{1}\frac{1}{m}\sum_{j=1}^{m}\frac
{t^{2j}(1+t)}{1-t}d\rho_{g}=0,
\]
and therefore condition (\ref{Mgen}) is satisfied. $\ \Diamond$

\bigskip

\textbf{Proof of Theorem \ref{shortmem}.}

\bigskip

Theorem \ref{shortmem} follows by combining Theorem \ref{short memory} with
the invariance principle in Kipnis and Varadhan (1997). We have only to verify
condition (\ref{key1}). It is known that the maximal inequality required by
condition (\ref{key1}) holds for partial sums of functions of reversible
Markov chains. Indeed, we know from Proposition 4 in Longla et al. (2012)
that
\begin{equation}
\mathbb{E(}\max_{1\leq i\leq n}S_{i}^{2})\leq2\mathbb{E(}\max_{1\leq i\leq
n}X_{i}^{2})+22\mathbb{\ }\max_{1\leq i\leq n}\mathbb{E(}S_{i}^{2})\label{LW}%
\end{equation}
and then, condition (\ref{condvar}) and stationarity implies condition
(\ref{key1}) with $c_{n}=\sqrt{n}$. $\ \Diamond$

\section{Appendix}

\textbf{Facts about spectral densities.} In the following lemma we combine a
few facts about spectral densities, covariances, behavior of variances of sums
and their relationships. The first two points are well known. They can be
found for instance in Bradley (2007, Vol 1, 0.19-0.21 and Ch.8). The point
(iii) was proven in Peligrad and Utev (2006-a).

\begin{lemma}
\label{spec} Let $(\xi_{i})_{i\in{\mathbb{Z}}}$ be a stationary sequence of
real valued variables with $\mathbb{E}\xi_{0}=0$ and finite second moment$.$
Let $F$ denotes the spectral measure and $f$ denotes its spectral density (if
exists) i.e.
\[
\mathbb{E}(\xi_{0}\xi_{k})=\int_{-\pi}^{\pi}e^{-ikt}dF(t)=\int_{-\pi}^{\pi
}e^{-ikt}f(t)dt.
\]
(i) For any positive integer $n$ and any real numbers $a_{1},\ldots,a_{n}$,
\begin{align*}
\mathbb{E}\left(  \sum_{k=1}^{n}a_{k}\xi_{k}\right)  ^{2} &  =\int_{-\pi}%
^{\pi}\left\vert \sum_{k=1}^{n}a_{k}e^{ikt}\right\vert ^{2}f(t)dt\leq2\pi\Vert
f\Vert_{\infty}\sum_{k=1}^{n}a_{k}^{2}\\
&  \leq\left(  \mathbb{E}[\xi_{0}^{2}]+2\sum_{k\geq1}|\mathbb{E}(\xi_{0}%
\xi_{k})|\right)  \sum_{k=1}^{n}a_{k}^{2}.
\end{align*}
(ii) Assume$\;\sum_{k=1}^{\infty}|\mathbb{E}(\xi_{0}\xi_{k})|<\infty$.\quad
Then, $f$ is continuous.\newline(iii) Assume that the spectral density $f$ is
continuous, and let $(d_{n,j})_{j\in{\mathbb{Z}}}$ be a double array of real
numbers with $d_{n}^{2}=\Sigma_{j\in{\mathbb{Z}}}d_{n,j}^{2}<\infty$ that
satisfies the condition
\begin{equation}
\frac{1}{d_{n}^{2}}\sum_{j\in\mathbb{Z}}|d_{n,j}-d_{n,j-1}|^{2}\rightarrow
0\text{.}\label{dcond}%
\end{equation}
Then,
\begin{equation}
\lim_{n\rightarrow\infty}\frac{1}{d_{n}^{2}}\mathbb{E}\left(  \sum
_{j\in\mathbb{Z}}d_{n,j}\xi_{j}\right)  ^{2}=2\pi f(0).\label{var}%
\end{equation}

\end{lemma}

\textbf{Acknowledgement.} The authors would like to thank the referee for
carefully reading the manuscript and for suggestions that improved the
presentation of this paper.

\end{document}